\documentclass{article}
\usepackage{amsmath,amsthm}

\usepackage{amssymb}

\newtheoremstyle{theorem}%name
  {10pt}          % space above
  {10pt}  % space below
  {\sl}  % bofy font
  {\parindent}     % ident - empty=no indent,  \parindent= paragraph indent
  {\bf}  % thm head font
  {. }    % punctuation after thm head
  { }    % space after thm head: `` ``=normal \newline=linebreak
  {}     % thm head specification
\theoremstyle{theorem}
\newtheorem{theorem}{Theorem}

\newtheorem{lemma}[theorem]{Lemma}

\newtheoremstyle{defi}%name
  {10pt}          % space above
  {10pt}  % space below
  {\rm}  % bofy font
  {\parindent}     % ident - empty=no indent,  \parindent= paragraph indent
  {\bf}  % thm head font
  {. }    % punctuation after thm head
  { }    % space after thm head: `` ``=normal \newline=linebreak
  {}     % thm head specification
\theoremstyle{defi}

\begin{document}

\title{Existence of Infinitely Many Solutions
for a Quasilinear Elliptic Problem on Time Scales}

\author{Moulay Rchid Sidi Ammi\\
\texttt{sidiammi@mat.ua.pt}
\and Delfim F. M. Torres\\
\texttt{delfim@mat.ua.pt}}

\date{Department of Mathematics\\
University of Aveiro\\
3810-193 Aveiro, Portugal}

\maketitle

%%%%%%%%%%%%%%%%%%%%%

\begin{abstract}
We study a boundary-value  quasilinear elliptic problem on a
generic time scale. Making use of the fixed-point index theory,
sufficient conditions are given to obtain existence, multiplicity,
and infinite solvability of positive solutions.

\smallskip

{\bf AMS Subject Classification:} 34B18, 39A10, 93C70.

{\bf Key Words and Phrases:} time scales, $p$-Laplacian, positive
solution, existence, multiplicity, infinite solvability.

\end{abstract}

%%%%%%%%%%%%%%%%%%%%%%%%%%%%%%%%%%%%%%%%%%%%%%%%%%%%%%%%%%%%%%%%%%%%

\section{Introduction}
\label{sec1}

We are interested in the study of the following quasilinear
elliptic problem:
\begin{equation} \label{eq1}
\begin{gathered}
-\left ( \phi_{p}(u^{\triangle}(t))\right)^{\nabla}= f(u(t)) +
h(t)\, , \quad  \forall t \in  (0,T)_{\mathbb{T}}= \mathbb{T} \, ,\\
u^{\triangle}(0)=0  \, ,\\
u(T) - u(\eta)=0 \, , \quad 0< \eta < T \, ,
\end{gathered}
\end{equation}
where $\phi_{p}(\cdot)$ is the $p-$Laplacian operator defined by
$\phi_{p}(s)= |s|^{p-2} s$, $p>1$, $(\phi_{p})^{-1}= \phi_{q}$
with $q$ the Holder conjugate of $p$, \textrm{i.e.} $\frac{1}{p}+
\frac{1}{q}= 1$, and $\mathbb{T}$ is a time-scale. We assume the
following hypotheses:
\begin{itemize}
\item[(H1)] The function $f: (0, T)\rightarrow \mathbb{R}^{+*}$ is
a continuous function;

\item[(H2)] The function $ h: T \rightarrow \mathbb{R}^{+*}$ is
left dense continuous (\textrm{i.e}, $h \in
\mathbb{C}_{ld}(\mathbb{T}, \mathbb{R}^{+*})$). Here
$\mathbb{C}_{ld}(\mathbb{T}, \mathbb{R}^{+*})$ denotes the set of
all left dense continuous functions from $\mathbb{T} \rightarrow
\mathbb{R}^{+}$; and $h \in L^{\infty}(0, T)$.
\end{itemize}

Results on existence of radially infinity many solutions for
\eqref{eq1} are proved using (i) variational methods \cite{ar,ck},
where solutions are obtained as critical points of some energy
functional on a Sobolev space, by imposing appropriate conditions
on $f$; (ii) methods based on phase-plane analysis and the
shooting method \cite{et}; (iii) by adapting the technique of time
maps \cite{eg}. For $p=2$, $h=0$, $\mathbb{T}=\mathbb{R}$, problem
\eqref{eq1} becomes a boundary-value problem of differential
equations. Our results extend and include results of the earlier
works to the case of a generic time-scale $\mathbb{T}$, $p\neq 2$
and where $h$ is not identically zero. In the case of $h=0$,
$p=2$, many existence results of dynamic equations on time scales
are available, using different fixed point theorems
\cite{and,kau}. We remark that there are not many results
concerning the $p-$Laplacian problems on time scales \cite{da}. In
this paper we prove existence of solutions by constructing an
operator whose fixed points are solutions of \eqref{eq1}. Our main
ingredient is the following well-known fixed-point theorem of
index theory.

\begin{theorem}[\cite{g1,gll}]
\label{thm11} Suppose  $E$ is a real Banach space, and $K\subset
E$ is a cone in $X$. Let  $\Omega _r= \left \{ u \in K, \|u\| < r
\right \}$, and $F: \Omega_r \to K$ be a completely continuous
operator satisfying $Fx \neq x$, for all $x \in \partial
\Omega_r$. The following holds:
\begin{enumerate}
\item[(i)] if $\|Fx\|\leq \|x\|$, $\forall x \in \partial \Omega
_r$, then $i(F, \Omega_r, K)=1$;

\item [(ii)] if $\|Fx\|\geq \|x\|$, $\forall x \in \partial
\Omega_r$, then $i(F, \Omega_r, K)=0$,
\end{enumerate}
where $i$ is the index of $F$.
\end{theorem}

The study of dynamic equations on time scales is a fairly new
subject, and research in this area is rapidly growing. The theory
of time scales has been created in order to unify continuous and
discrete analysis, allowing a simultaneous treatment of
differential and difference equations, and to extend those
theories to so-called delta/nabla-dynamic equations. A vast
literature has already emerged in this field: see \textrm{e.g.}
\cite{abra,a3,a4}. For an introduction to time scales with
applications, we refer the reader to \cite{bp1,bp3}.

The outline of the paper is as follows. In Section~\ref{sec2} we
give some preliminary results with respect to the calculus on time
scales.  Section~\ref{sec3} is devoted to the existence of
positive solutions using fixed-point index theory. The remaining
sections deal with multiplicity and infinite solvability solutions
for \eqref{eq1}.

%%%%%%%%%%%%%%%%%%%%%%

\section{Preliminary results on time scales}
\label{sec2}

We begin by recalling some basic concepts of time scales. Then, we
prove some preliminary results that will be needed in the sequel.

A time scale $\mathbb{T}$ is an arbitrary nonempty closed subset
of the set  $\mathbb{R}$ of real numbers. The operators $\sigma$
and $\rho$ from $\mathbb{T}$ to $\mathbb{T}$ are defined in
\cite{h1,h2}:
$$
\sigma(t)=\inf\{\tau\in\mathbb{T}:\tau> t\}\in\mathbb{T}, \quad
\rho(t)=\sup\{\tau\in\mathbb{T}:\tau< t\}\in\mathbb{T}.
$$
They are called the  forward jump operator and the backward jump
operator, respectively.

The point $t\in\mathbb{T}$ is said to be left-dense,
left-scattered, right-dense, or right-scattered, if $\rho (t)=t$,
$\rho(t)<t$, $\sigma(t)=t$, $\sigma(t)>t$, respectively. If
$\mathbb{T}$ has a right scattered minimum $m$, we define
$\mathbb{T}_{k}=\mathbb{T}-\{m\}$; otherwise we set
$\mathbb{T}_{k}=\mathbb{T}$. Similarly, if $\mathbb{T}$ has a left
scattered maximum $M$, we define
$\mathbb{T}^{k}=\mathbb{T}-\{M\}$; otherwise we set
$\mathbb{T}^{k}=\mathbb{T}$.

Let $f:\mathbb{T} \to \mathbb{R}$ and $t\in \mathbb{T}^{k}$
(assume $t$ is not left-scattered if $t=\sup\mathbb{T}$). We
define $f^{\triangle}(t)$ to be the number (provided it exists)
such that given any $\epsilon>0$ there is a neighborhood $U$ of
$t$ such that
$$
|f(\sigma(t))-f(s)-f^{\Delta}(t)(\sigma(t)-s) |\le | \sigma(t)-s|,
\quad \mbox{for all } s\in U.
$$
We call $f^{\triangle}(t)$ the delta derivative of $f$ at $t$. We
remark that $f^{\triangle}$ is the usual derivative $f'$ if
$\mathbb{T}=\mathbb{R}$ and the usual forward difference $\Delta
f$ (defined by $\triangle f(t)=f(t+1)-f(t))$ if
$\mathbb{T}=\mathbb{Z}$.

Similarly, for $t\in \mathbb{T}$ (assume $t$ is not
right-scattered if $t=\inf\mathbb{T}$),  the nabla derivative of
$f$ at the point $t$ is defined in \cite{a4} to be the number
$f^{\nabla}(t)$ (provided it exists) with the property that for
each $\epsilon >0$ there is a neighborhood $U$ of $t$ such that
$$
|f(\rho(t))-f(s)-f^{\nabla}(t)(\rho(t)-s) |\le | \rho(t)-s |, \quad
\mbox{for all } s\in U.
$$
If $\mathbb{T}=\mathbb{R}$, then $f^\Delta(t)=f^\nabla(t)=f'(t)$. If
$\mathbb{T}=Z$, then $f^\nabla(t)=f(t)-f(t-1)$ is the backward
difference operator.

We say that a function $f$ is left-dense continuous
($ld$-continuous for short), provided $f$ is continuous at each
left-dense point in $\mathbb{T}$ and its right-sided limit exists
at each right-dense point in $\mathbb{T}$. It is well-known that
if $f$ is $ld$-continuous and if $F^{\nabla}(t)=f(t)$, then one
can define the nabla integral by
$$
\int^b_a f(t)\nabla t=F(b)-F(a).
$$
If $F^{\Delta}(t)=f(t)$, then we define the delta integral by
$$
\int^b_a f(t)\Delta t=F(b)-F(a).
$$
Function $F$ is said to be an antiderivative of $f$. For more
details on time scales, the reader can consult
\cite{a1,abra,a2,a3,bp1,bp3} and references therein.

In the rest of the paper, $\mathbb{T}$ is a closed subset of
$\mathbb{R}$ with $0\in\mathbb{T}_k$, $T\in\mathbb{T}^k$. We
denote $E= \mathbb{C}_{ld}([0, T], \mathbb{R})$, which is a Banach
space with the maximum norm $\|u\|= \max_{[0, T]}|u(t)|$.

\begin{lemma}
\label{lm1} Suppose that conditions $(H1)$ and $(H2)$ hold. Then,
$u(t)$ is a solution of the boundary-value problem \eqref{eq1} if
and only if $u(t) \in E$ is a solution of the following equation:
\begin{equation}\label{eq2}
 u(t)= \phi_{q} \left (\int_{\eta}^{T}(f(u(r)+h(r)) \, \nabla r
 \right )
+\int_{0}^{t}\phi_{q}\left ( \int_{s}^{T} (f(u(r)+h(r)) \, \nabla r
\right ) \triangle s.
\end{equation}
\end{lemma}

\begin{proof}
By integrating the equation \eqref{eq1} on $(s, T)$, we have
$$
\phi_{p}(u^{\triangle}(T))= \phi_{p}(u^{\triangle}(s))-
\int_{s}^{T}(f(u(r)+h(r)) \, \nabla r.
$$
Then,
$$
\phi_{p}(u^{\triangle}(s))= \phi_{p}(u^{\triangle}(T))+
\int_{s}^{T}(f(u(r)+h(r)) \, \nabla r.
$$
Using the boundary conditions, we have
$$
\phi_{p}(u^{\triangle}(s))= \int_{s}^{T}(f(u(r)+h(r)) \, \nabla r.
$$
Thus,
\begin{equation*}
u^{\triangle}(s)= \phi_{q} \left ( \int_{s}^{T}(f(u(r)+h(r)) \,
\nabla r \right ).
\end{equation*}
Integrating the last equation on
$(0, t)$, we have
\begin{equation*}
\begin{split}
u(t)&=u(0) + \int_{0}^{t} \phi_{q} \left ( \int_{s}^{T}(f(u(r)+h(r))
\, \nabla r \right ) \triangle s\\
& = u^{\triangle}(\eta) + \int_{0}^{t} \phi_{q} \left (
\int_{s}^{T}(f(u(r)+h(r))
\, \nabla r \right ) \triangle s\\
&= \phi_{q} \left ( \int_{\eta}^{T} (f(u(r)+h(r)) \, \nabla r \right
) + \int_{0}^{t} \phi_{q} \left ( \int_{s}^{T}(f(u(r)+h(r)) \,
\nabla r \right ) \triangle s.
\end{split}
\end{equation*}
Inversely, if we suppose that \eqref{eq2} holds, it is easy to get
the first equation of \eqref{eq1} by derivation, and to see that
$u$ satisfies the boundary value  conditions in \eqref{eq1}.
Furthermore, $u$ is obviously positive since $\phi_{q}$  is non
decreasing function and $f$ and $h$ are also  positives functions.
The proof of Lemma~\ref{lm1} is now complete.
\end{proof}

On the other hand, we have $-(\phi_{p}(u^{\triangle}))^{\nabla} =
f(u(t))+ h(t)$. Then, since $f, h \geq 0$, we have
$(\phi_{p}(u^{\triangle}))^{\nabla} \leq 0$ and
$(\phi_{p}(u^{\triangle}(t_2))) \leq
(\phi_{p}(u^{\triangle}(t_1)))$ for any $t_1,\, t_2 \in [0, T]$
with $t_1 \leq t_2$. It follows that $u^{\triangle}(t_2) \leq
u^{\triangle}(t_1)$ for $t_1 \leq t_2$. Hence, $u^{\triangle} (t)$
is a decreasing function on $[0, T]$. This means that the graph of
$u(t)$ is concave on $[0, T]$.

Let $K \subset E$ be the cone defined by
$$
K= \left \{ u \in E: u(t)\geq 0, u(t) \mbox{ is a concave
function}, t\in [0, 1] \right \},
$$
and $F: K\rightarrow E$ the operator
$$
Fu(t) = \phi_{q} \left ( \int_{\eta}^{T} (f(u(r))+h(r)) \, \nabla r
\right )+ \int_{0}^{t} \phi_{q} \left ( \int_{s}^{T}(f(u(r)+h(r)) \,
\nabla r \right ) \triangle s.
$$
It is easy to see that \eqref{eq1} has a solution $u= u(t)$ if and
only if $u$ is a fixed point of the operator $F$. One can also
verify that $F(K) \subset K$ and $F: K \rightarrow K$ is
completely continuous.

%%%%%%%%%%%%%%%%%%%%%%%%%%%%%%%%%%%%%%%%%%%%%%%%

\section{ Existence of positive solutions}
\label{sec3}

We define two open subsets $\Omega_1$ and $\Omega_2$ of $E$:
$$
\Omega_1 = \left \{ u \in K: \|u\| < a \right\} \, , \quad
\Omega_2 = \left \{ u \in K: \|u\| < b \right\}.
$$
Without loss of generality, we suppose that $b < a$. For
convenience, we introduce the following notation:
\begin{equation*}
\begin{split}
A& = \frac{a- \alpha \|h\|_{\infty}^{1/p-1}}{\alpha a},\mbox{ where
}
\alpha = \phi_{q}(2^{p-2}) \phi_{q}(T)(T+1),\\
B& = \phi_{p}(T- \eta).
\end{split}
\end{equation*}

\begin{theorem}\label{thm31}
Besides $(H1)$ and $(H2)$, suppose that $f$ also satisfies:

\begin{description}

\item[(i)] $\max_{0\leq u \leq a}f(u) \leq \phi_{p}(aA)$;

\item[(ii)] $\min_{0\leq u \leq b}f(u) \geq
\phi_{p}(bB)$.

\end{description}
Then, \eqref{eq1} has a positive solution.
\end{theorem}

\begin{proof}
If $u \in \partial \Omega_{1}$, we have:
\begin{equation*}
\begin{gathered}
\begin{split}
\|F(u)\| &  \leq \phi_{q} \left ( \int_{\eta}^{T}((aA)^{p-1}+
\|h\|_{\infty}) \nabla r \right ) \\
& \qquad+ \int_{0}^{T} \phi_{q} \left
(\int_{s}^{T}((aA)^{p-1}+ \|h\|_{\infty}) \nabla r \right )
\triangle s
\end{split}\\
\begin{split}
&\leq \phi_{q}\left ( ((aA)^{p-1}+ \|h\|_{\infty})(T- \eta) \right) \\
& \qquad + \int_{0}^{T} \phi_{q} \left ( ((aA)^{p-1}+ \|h\|_{\infty})(T- s)
\right ) \triangle s
\end{split}
\end{gathered}
\end{equation*}
\begin{equation*}
\begin{gathered}
\begin{split}
&\leq \phi_{q}\left ( (aA)^{p-1}+ \|h\|_{\infty}\right )
\phi_{q}(T) \\
& \qquad + \phi_{q}\left ( (aA)^{p-1}+ \|h\|_{\infty}\right )
\int_{0}^{T} \phi_{q}(T-s) \triangle s
\end{split}\\
\begin{split}
&\leq \phi_{q}\left ( (aA)^{p-1}+
(\|h\|_{\infty}^{1/p-1})^{p-1}\right ) \phi_{q}(T) \\
& \qquad + \phi_{q}\left (
(aA)^{p-1}+ (\|h\|_{\infty}^{1/p-1})^{p-1}\right )
\phi_{q}(T) T
\end{split}\\
\leq \phi_{q}\left ( (aA)^{p-1}+
(\|h\|_{\infty}^{1/p-1})^{p-1}\right ) \phi_{q}(T) (T+1).
\end{gathered}
\end{equation*}
Using the elementary inequality
$$
x^{p} + y^{p} \leq 2^{p-1}(x+y)^{p},
$$
and the form of $A$, it follows that
\begin{equation*}
\begin{split}
\|F(u)\| &  \leq \phi_{q}(T+1)(2^{p-2})(aA+
\|h\|_{\infty}^{1/p-1})\\
& \leq \|u\| = a.
\end{split}
\end{equation*}
Therefore, $\|Fu\| \leq \|u\|$ for all $u \in \partial
\Omega_{1}$. Then, by Theorem~\ref{thm11},
\begin{equation}
\label{eq4}
i(F, \Omega_{1}, K)= 1.
\end{equation}
On the other hand, for $u \in \partial \Omega_{2}$ we have:
\begin{equation*}
\begin{split}
\|F(u)\| & \geq \phi_{q}\left ( \int_{\eta}^{T} f(u(r) \, \nabla r
\right ) + \int_{0}^{t} \phi_{q} \left ( \int_{s}^{T}(f(u(r)) \,
\nabla r
\right ) \triangle s. \\
& \geq B b \phi_{q} (T- \eta ) \\
&\geq b=\|u\| (\mbox{ since } B= \phi_{p}(T- \eta )).
\end{split}
\end{equation*}
Therefore, $\|Fu\| \geq \|u\|$ for all $u \in \partial
\Omega_{2}$. By Theorem~\ref{thm11},
\begin{equation}\label{eq3}
i(F, \Omega_{2}, K)= 0.
\end{equation}
It follows by \eqref{eq4} and \eqref{eq3} that
$$
i(F, \Omega_{1}\backslash \overline{\Omega_{2}}, K)= 1.
$$
Then $T$ has a fixed point $u \in \Omega_{1}\backslash
\overline{\Omega_{2}}$. Obviously, $u$ is a positive solution of
problem \eqref{eq1} and $b < \|u\| < a$. The proof of
Theorem~\ref{thm31} is complete.
\end{proof}

%%%%%%%%%%%%%%%%%%%%%%%%%%%%%%%%%%%%%%%%%%%%%%%%%%%%%%%%%%%%

\section{Multiplicity}
\label{sec4}

By multiplicity we mean the existence of an arbitrary number of
solutions. We now obtain results on the multiplicity of positive
solutions for \eqref{eq1} under the following assumptions: we
suppose that there exist positive real numbers $0 <a_{1}<a_{2}<
\ldots < a_{k+1}$, such that
\begin{description}

\item[$(i)$] $\max_{0\leq u \leq a_{2i-1}}f(u) \leq
\phi_{p}(a_{2i-1}A), i= 1, \ldots [\frac{k+2}{2}]$;

\item[$(ii)$] $\min_{0\leq u \leq a_{2i}}f(u) \geq
\phi_{p}(a_{2i}B), i= 1, \ldots [\frac{k+1}{2}]$;

\end{description}
where $[n]$ denote the integer part of $n$.

\begin{theorem} Assume that $(i)$-$(ii)$ hold.
Then, problem \eqref{eq1} has at least $k$ positive solutions
$u_{1}, \ldots u_{k}$ such that
$$
a_{i} < \|u_{i}\| < a_{i+1}, \quad i=1, \ldots k  \, .
$$
\end{theorem}

\begin{proof}
By continuity, there exist
$$
0< b_1 < a_1 < c_1 < b_2 < a_2 < c_2< \ldots c_k
 < b_{k+1} < a_{k+1} < + \infty
$$
such that
$$
\min_{0 \leq u \leq b_{2i-1}} f(u) \geq \phi_{p}(b_{2i-1 }B) \, ,
\quad \min_{0 \leq u \leq c_{2i-1}} f(u) \geq \phi_{p}(c_{2i-1
}B),
$$
for $ i= 1, \ldots [\frac{k+2}{2}]$, and
$$
\max_{0 \leq u \leq c_{2i}} f(u) \leq \phi_{p}(c_{2i}A) \, , \quad
\max_{0 \leq u \leq b_{2i}} f(u) \leq \phi_{p}(b_{2i }A),
$$
for $i= 1, \ldots [\frac{k+1}{2}]$. Then, calling
Theorem~\ref{thm31} to each interval $(c_{i}, b_{i+1}), i= 1,
\ldots k$, we obtain the intended result.
\end{proof}

%%%%%%%%%%%%%%%%%%%%%%%%%%%%%%%%%%%%%%%%%%%%%%%%%%%%%%%%%%%%

\section{Infinite solvability}
\label{sec5}

\begin{theorem}
If the following two conditions hold,
\begin{description}

\item[$(i)$] $ \liminf_{a \rightarrow 0} \frac{ \max_{0\leq u \leq
a} \{f(u)\}}{a^{p-1}} \leq \phi_{p}(A)$,

\item[$(ii)$] $ \limsup_{b \rightarrow 0} \frac{ \max_{0\leq u
\leq b} \{f(u) \}}{b^{p-1}} \geq \phi_{p}(B)$,
\end{description}
then, problem \eqref{eq1} has a sequence of positive solutions
$(u_{k})_{k\geq 1}$ such that $\|u_{k}\| \rightarrow 0$ as $k
\rightarrow \infty$.
\end{theorem}

\begin{proof}
From $(i)$ and $(ii)$, there exists a sequence of pairs of
positive numbers $(a_{k}, b_{k})_{k \geq 1}$ convergent to $(0,
0)$ such that
$$
\max_{0 \leq u \leq a_{k}} f(u) \leq \phi_{p}(a_{k}A),
$$
$$
\min_{0 \leq u \leq b_{k}} f(u) \geq \phi_{p}(b_{k}B).
$$
Suppose that
$$
a_{1} > b_{1} >a_{2} >b_{2} > \ldots a_{k} >b_{k} > \ldots
$$
Calling Theorem~\ref{thm31} on each pair $(a_{k}, b_{k})_{k \geq
1}$, we conclude that \eqref{eq1} has a sequence of positive
solutions $(u_{k})_{k\geq 1}$ such that $b_{k} \leq \|u_{k}\| \leq
a_{k}$.
\end{proof}

%%%%%%%%%%%%%%%%%%%%%

\section*{Acknowledgements}

The authors were supported by the \emph{Portuguese Foundation for
Science and Technology} (FCT) through the \emph{Centre for Research
in Optimization and Control} (CEOC) of the University of Aveiro,
cofinanced by the European Community fund FEDER/POCTI, and the
project SFRH/BPD/20934/2004.

%%%%%%%%%%%%%%%%%%%%%

%%%%%%%%%%%%%%%%%%%%%

\end{document}